\documentclass{amsproc}

\usepackage{amssymb}
\usepackage{graphicx}
\usepackage{amscd}
\usepackage{amsmath}
\theoremstyle{plain}

\newtheorem{lemma}{Lemma}

\newtheorem{proposition}{Proposition}

\numberwithin{equation}{section}

\begin{document}

\title[On the dynamic representation of Weyl solution]{Dynamic representation of the Weyl solution for the Schr\"odinger operator on the semi-axis}

\author{A. S. Mikhaylov}
\address{St. Petersburg   Department   of   V.A. Steklov    Institute   of   Mathematics
of   the   Russian   Academy   of   Sciences, 7, Fontanka, 191023
St. Petersburg, Russia and Saint Petersburg State University,
St.Petersburg State University, 7/9 Universitetskaya nab., St.
Petersburg, 199034 Russia.}
\email{mikhaylov@pdmi.ras.ru}

\author{ V. S. Mikhaylov}
\address{St.Petersburg   Department   of   V.A.Steklov    Institute   of   Mathematics
of   the   Russian   Academy   of   Sciences, 7, Fontanka, 191023
St. Petersburg, Russia.} \email{ftvsm78@gmail.com}

\thanks{No new data were created or analyzed during this study. Data sharing is not applicable to this article.}

\thanks{All authors declare that they have no conflicts of interest.}

\date{November, 2025}

\maketitle






\noindent {\bf Abstract.} We derive a representation formula for the Weyl solution
to the Schr\"odinger operator on the semi-axis for certain classes of potentials.
Our approach is based on relations with the initial-boundary value problem
for the wave equation with the same potential on the half-line.

\section{Introduction.}

We consider the Schr\"odinger operator
\begin{equation}
H=-\partial _{x}^{2}+q(x),\quad x>0  \label{Schrodinger_op}
\end{equation}
on $L_{2}(\mathbb{R}_{+})$ with a real-valued locally integrable potential $q.$

We assume that (\ref{Schrodinger_op}) is limit point case at $\infty $, that is, for each $z\in
\mathbb{C}_{+}:=\{z\in \mathbb{C}:\operatorname{Im}z>0\}$ the equation
\begin{equation}
\label{Schrod_eqn}
-u^{\prime \prime }+q\left( x\right) u=zu
\end{equation}
has a unique solution $u_{+}$ which is in $L_{2}$:
\begin{equation*}
\int_{\mathbb{R}_{+}}\left\vert u_{+}\left( x,z\right) \right\vert
^{2}dx<\infty ,z\in \mathbb{C}_{+}.
\end{equation*}
Such a solution $u_{+}$ is called a \emph{Weyl} \emph{solution} and its
existence is the central point of the Titchmarsh-Weyl theory.

The \emph{Titchmarsh-Weyl }$m$\emph{-function}, $m(z),$ is defined for $z\in \mathbb{C}_{+}$ as
\begin{equation}
\label{Weyl_def}
m\left( z\right) :=\frac{u_{+}^{\prime }\left( 0,z\right) }{u_{+}\left(
0,z\right) }.
\end{equation}
The function $m(z)$ is analytic in $\mathbb{C}_{+}$ and satisfy the Herglotz property, i.e. $m: \mathbb{C}_+\mapsto \mathbb{C}_+$.

Simon and Gesztey in \cite{S99,GS00} proposed the following representation for the Weyl function:
\begin{equation}
\label{Simon_rep}
m(-k^{2})=-k-\int_{0}^{\infty }A(\alpha )e^{-2\alpha k}\,d\alpha,
\end{equation}
where the absolute convergence of integral was proven for $q\in
L^{1}\left( \mathbb{R}_{+}\right) $ and $q\in L^{\infty }\left( \mathbb{R}%
_{+}\right) $. They called the function $A$ in (\ref{Simon_rep}) the\emph{\ }$A-$\emph{amplitude}.

In \cite{AMR} the authors proposed the natural physical interpretation of the $A-$amplitude.
It is as follows: for the same potential $q$, the initial-value problem for the wave equation is considered:
\begin{equation}
\label{wave_eqn}
\begin{cases}
u_{tt}(x,t)-u_{xx}(x,t)+q(x)u(x,t)=0, \quad x>0,\ t>0,\\
u(x,0)=u_t(x,0)=0,\ u(0,t)=f(t).
\end{cases}
\end{equation}
where $f$ is an arbitrary $L^2_{loc}\left( \mathbb{R}_{+}\right) $
function referred to as a \emph{boundary control}. The solution to (\ref{wave_eqn}) is denoted by $u^f$.
The {\bf response operator} (the dynamical Dirichlet-to-Neumann
map) $R^T$ for the system (\ref{wave_eqn}) is defined in
$\mathcal{F}^T:=L_2(0,T)$ by
\begin{equation*}
(R^Tf)(t)=u_{x}^f(0,t), \ t \in (0,T),
\end{equation*}
with the domain $\{ f \in C^2([0,T]):\; f(0)=f'(0)=0\}.$ The
response operator has the representation (see \cite{AM,AMR}):
\begin{equation*}
(R^Tf)(t)=-f'(t)+\int_0^tr(s)f(t-s)\,ds,
\end{equation*}
where $r$ is called the \emph{response function}. In \cite{AMR} the authors have shown that
\begin{eqnarray}
m(-k^2)=-k+\int_0^\infty r(\alpha)e^{-k\alpha},\label{Resp_rep}\\
A(\alpha)=-2r(2\alpha) \label{A_and_r}.
\end{eqnarray}
and the integrals in (\ref{Simon_rep}) and (\ref{Resp_rep})
converges for wider class of potentials, namely for $q\in l^\infty(L_1(\mathbb{R}_+)):=\left\{q\,|\, \int_n^{n+1}|q(x)|\,dx\in l^\infty\right\}$.
Formula (\ref{A_and_r}) shows that $A-$amplitude, initially introduced in some artificial way, is in fact (taking into account the sign and scaling of coordinates) a response function,
i.e., the kernel of the response operator, a classical object in the theory of dynamic inverse problems \cite{B07}.

Although the authors in \cite{AMR} pointed out the connection between (\ref{Schrod_eqn}) and (\ref{wave_eqn}) using the
Laplace transform, they dealt directly with the integral equation for the $A-$amplitude and did not show that the solution to (\ref{wave_eqn}) after applying the
Laplace transform is in fact the Weyl solution. This drawback can be overcome by appealing to the finite wave propagation speed in the system (\ref{wave_eqn}),
as was done for the system associated with Jacobi matrices in \cite{MMS,MM5} (note that, nevertheless in \cite{MM6} the authors proved that the corresponding solution is from $l_2$),
but here we propose a different
approach and show how to obtain the Weyl solution for (\ref{Schrod_eqn}) from the solution $u^f$ of (\ref{wave_eqn}). Thus, the goal of this paper is to obtain
a new representation for the Weyl solution which is interesting in itself and eliminates the drawback from \cite{AMR}. Note that our results do not follow from \cite{AMR},
but rather the opposite.

The paper is organized as follows: in the second section we follow
\cite{AM}: formulate "standard" statements on the solvability of
(\ref{wave_eqn}) depending on the potential, show a representation
of $u^f$ in terms of the solution of a certain Goursat problem. We
rewrite this Goursat problem as an integral equation and formulate
some statements about its solvability for a potential from
$L_{1,\operatorname{loc}(\mathbb{R}_+)}$. In the third section, we
show that for $L_1(\mathbb{R}_+)$ potentials, a representation for
the Weyl solution follows from the "standard" estimates from
\cite{AM}. We then analyze the integral equation for the solution
of the Goursat problem in the case of $q\in
l^\infty\left(L_1(\mathbb{R}_+)\right)$, which also provides a
representation for the Weyl solution in this case.

\section{Dynamical system.}

In \cite{AM} the authors shown that the solution to (\ref{wave_eqn}) admits the following representation
\begin{equation}
\label{wave_eqn_sol} u^f(x,t)=\begin{cases}
f(t-x)+\int_x^tw(x,s)f(t-s)\,ds, \quad x \leq t,\\
0, \quad x > t.\end{cases}.
\end{equation}
Where the kernel $w(w,t)$ is the unique solution to
the Goursat problem:
\begin{equation}
\label{gursa}
\begin{cases}
w_{tt}(x,t)-w_{xx}(x,t)+q(x)w(x,t)=0, \quad 0<x<t,\\
w(0,t)=0,\ w(x,x)=-1/2\int_0^xq(s)\,ds.
\end{cases}
\end{equation}
For $u^f$ they proved the following
\begin{proposition}
\label{Prop_wave_sol}
\begin{itemize}
\item[a)] If $q\in C^1(\mathbb{R}_+)$, $f\in C^2(\mathbb{R}_+)$
and $f(0)=f'(0)=0$, then the solution to (\ref{wave_eqn}) given by formula (\ref{wave_eqn_sol})
is a classical solution to (\ref{wave_eqn}).
\item[b)] If $q\in
L^1_{loc}(\mathbb{R}_+)$ and $f\in L^2(0,T)$, then formula
$(\ref{wave_eqn_sol})$ represents a unique generalized solution to
the initial-boundary value problem $(\ref{wave_eqn})$  \\ $u^f \in
C([0,T];\mathcal{H}^T),$ where  $\mathcal{H}=L^2_{loc}(0,\infty)$ \ \
\mbox{and} \ \ $\mathcal{H}^T := \{u \in \mathcal{H}: \, \mbox{
supp } u \subset [0,T]\, \}.$
\end{itemize}
\end{proposition}

\begin{proposition}
\label{Prop_Goursat}
\begin{itemize}
\item[a)] If $q\in L^1_{loc}(\mathbb{R}_+)$, then the generalized
solution $w(x,s)$ to the Goursat problem $(\ref{gursa})$ is a
continuous function and
\begin{eqnarray}
|w(x,s)|\leqslant
\Bigl(\frac{1}{2}\int_0^{\frac{s+x}{2}}|q(\alpha)|\,d\alpha\Bigr)
\operatorname{exp}\Bigl\{\frac{s-x}{4}\int_0^{\frac{s+x}{2}}|q(\alpha)|\,d\alpha\Bigr\}\label{gursa_est},\\
w_x(\cdot,s), w_s(\cdot,s), w_x(x,\cdot), w_s(x,\cdot)\in
L_{1,\,loc}(\mathbb{R}_+).\label{gursa_deriv}
\end{eqnarray}
Partial derivatives in $(\ref{gursa_deriv})$ continuously in
$L^1_{loc}(\mathbb{R}_+)$ depend on parameters $x$, $s$. The
equation in $(\ref{gursa})$ holds almost everywhere and the
boundary conditions are satisfied in the classical sense.

\item[b)] If $q\in C_{loc}(\mathbb{R}_+)$, then the generalized
solution $ w(x,s)$ to the Goursat problem $(\ref{gursa})$ is
$C^1$-smooth, equation and boundary conditions are satisfied in
the classical sense.

\item[c)] If $q\in C^1_{loc}(\mathbb{R}_+)$, then the solution to
the Goursat problem $(\ref{gursa})$ is classical, all its
derivatives up to the second order are continuous.
\end{itemize}
\end{proposition}

We outline the scheme of the proof of the Proposition \ref{Prop_Goursat}, which is given in \cite{AM}

By setting $\xi=s-x$, $\eta=s+x$, and
\begin{equation*}
v(\xi,\eta)=w\Bigl(\frac{\eta-\xi}{2},\frac{\eta+\xi}{2}\Bigl),
\end{equation*}
equation $(\ref{gursa})$ reduces to
\begin{equation}
\label{gursa_trans} \left\{
\begin{array}l
v_{\xi\eta}-\frac{1}{4}q(\frac{\eta-\xi}{2})v=0, \quad 0<\xi<\eta,\\
v(\eta,\eta)=0,\
v(0,\eta)=-\frac{1}{2}\int_0^{\eta/2}q(\alpha)\,d\alpha.
\end{array}
\right.
\end{equation}
The boundary value problem $(\ref{gursa_trans})$ is equivalent to the
integral equation
\begin{equation}
\label{wave_gursa_new_2}
v(\xi,\eta)=-\frac{1}{2}\int_{\xi/2}^{\eta/2}q(\alpha)\,d\alpha-\frac{1}{4}
\int_0^\xi\,d\xi_1\int_\xi^\eta\,d\eta_1q\Bigl(\frac{\eta_1-\xi_1}{2}\Bigr)v(\xi_1,\eta_1).
\end{equation}
We introduce a new function
\begin{equation}
\label{Q_func}
Q(\xi,\eta)=-\frac{1}{2}\int_{\xi/2}^{\eta/2}q(\alpha)\,d\alpha
\end{equation}
and the operator $K:C(\mathbb{R}^2)\mapsto C(\mathbb{R}^2)$ by the
rule
\begin{equation*}
\Bigl(Kv\Bigr)(\xi,\eta)=\frac{1}{4}
\int_0^\xi\,d\xi_1\int_\xi^\eta\,d\eta_1q\Bigl(\frac{\eta_1-\xi_1}{2}\Bigr)v(\xi_1,\eta_1).
\end{equation*}
Rewriting $(\ref{wave_gursa_new_2})$ as
\begin{equation*}
v=Q-Kv
\end{equation*}
and formally solving it by iterations, we get
\begin{equation}
\label{gursa_sol_sum} v(\xi,\eta)=Q(\xi,\eta)+\sum_{n=1}^\infty
(-1)^n(K^nQ)(\xi,\eta).
\end{equation}
To prove the convergence of $(\ref{gursa_sol_sum})$ we need
suitable estimates for $|K^nQ|(\xi,\eta)$. In \cite{AM} the
authors obtain such estimates for the case of $q\in
L_{1\operatorname{loc}(\mathbb{R}_+)}$, in the third Section we
refine these estimates for $q\in l^\infty(L_1(\mathbb{R}_+))$.

\section{Weyl solution representation. }

For the special control $f(t)=\delta(t)$, we denote by $u^\delta(x,t)$ the fundamental solution of (\ref{wave_eqn}). From (\ref{wave_eqn_sol})
it is clear that any solution of (\ref{wave_eqn}) has the form $u^f=u^\delta*f$.

Let $f \in C_0^{\infty}(0,\infty)$ and
\begin{equation*}
\widehat{f}(k):=\int_0^{\infty} f(t)\,e^{ikt}\,dt\,
\end{equation*}
be its Fourier transform. Function $\widehat{f}(k)$ is well
defined for $k \in \mathbb{C}$.
Going formally in (\ref{wave_eqn}) over to the Fourier transforms
(we note that the authors in \cite{S99,GS00,AMR} used the Laplace
transform), one has
\begin{equation*}
\begin{cases}
-\widehat{u^f_{xx}}(x,k)+q(x)\widehat{u^f}(x,k) =k^{2}\widehat{u^f}(x,k),\\
\widehat{u^f}(0,k) =\widehat{f}(k),
\end{cases}
\end{equation*}
and for the response operator
\begin{equation*}
\widehat{(Rf)}(k)=\widehat{u^f_x}(0,k)\,,
\end{equation*}
respectively.
We note that in view of (\ref{wave_eqn_sol}) $\widehat{u^f}(x,k)=\widehat{u^\delta}(x,k)\widehat f(k)$, so we see that the
values of the function $\widehat{u} (0,k) $ and its first derivative at the
origin, $\widehat{u}_x (0,k), $ are related through the Titchmarsh-Weyl
m-function (cf. (\ref{Weyl_def})):
\begin{equation*}
\widehat{u^f}_x (0,k)=m(k^2)\widehat{f}(k)\,.
\end{equation*}
Thus the formula (\ref{Resp_rep}) is valid provided we can justify
two things: taking the Fourier transform of (\ref{wave_eqn}) and
that $\widehat {u^{\delta}}(x,k)$ is a Weyl solution, where
\begin{equation}
\label{Weyl_func}
\widehat {u^{\delta}}(x,k)=e^{ikx}+\int_0^\infty w(x,t)e^{ikt}\,dt.
\end{equation}

\subsection{Case of $L_1(\mathbb{R}_+)$ potentials }

We note that the first term in (\ref{Weyl_func}) is in $L_2(\mathbb{R}_+)$ as soon as $\operatorname{Im}k>0$. Therefore, we need to estimate the second term in (\ref{Weyl_func}).
Using (\ref{gursa_est}) we obtain:
\begin{equation}
\label{W_est1}
|w(x,t)|\leqslant \frac{1}{2}\|q\|_{L_1}e^{\frac{\|q\|_{L_1}}{4}(t-x)}.
\end{equation}
The last estimate implies that $\int_0^\infty w(x,t)e^{ikt}\,dt$ converges as $
\operatorname{Im}k>\frac{\|q\|_{L_1}}{4}$.
Then from the same estimate (\ref{W_est1}) it follows that the second term in (\ref{Weyl_func}) is from $L_2(\mathbb{R}_+)$. Summarizing all of the above we obtain
\begin{proposition}
If $q\in L_1(\mathbb{R}_+)$ then formula (\ref{Weyl_func}) yields
the Weyl solution in the domain
$\operatorname{Im}k>\frac{\|q\|_{L_1}}{4}$.
\end{proposition}

\subsection{Case of $l_\infty(l_1(\mathbb{R}_+))$ potentials}

The following Lemma was proved in \cite{AMR}:
\begin{lemma}
\label{L1}
Let $f\left( x\right) $ be a non-negative function and
\begin{equation*}
||f||:=\sup_{x\geq 0}\int_{x}^{x+1}f\left( s\right) ds<\infty .  
\end{equation*}
Then for any $a,b\geq 0$ and natural $n$%
\begin{equation*}
\int_{0}^{a}\left( x+b\right) ^{n}f\left( x\right) dx\leq \frac{\left(
a+b+1\right) ^{n+1}}{n+1}||f||.  
\end{equation*}
\end{lemma}

We need to show the convergence of (\ref{gursa_sol_sum}), to this aim we rewrite
\begin{equation*}
Q(\xi,\eta)=-\frac{1}{2}\int_{\xi/2}^{\eta/2}q(\alpha)\,d\alpha=\left[\alpha=\frac{\gamma}{2}\right]=-\frac{1}{4}\int_\xi^\eta q\left(\frac{\gamma}{2}\right)\,d\gamma.
\end{equation*}
Then we introduce the notation $\widetilde q(\gamma):=q\left(\frac{\gamma}{2}\right)$ and use Lemma \ref{L1} to estimate:
\begin{equation*}
|Q(\xi,\eta)|\leqslant\frac{1}{4}\int_\xi^\eta \left|q\left(\frac{\gamma}{2}\right)\right|\,d\gamma\leqslant \frac{1}{4}(\eta+1)\|\widetilde q\|.
\end{equation*}
The second term in (\ref{gursa_sol_sum}) is estimated by:
\begin{eqnarray*}
|KQ(\xi,\eta)|\leqslant\frac{1}{4}\int_0^\xi\, d\xi_1 \int_\xi^\eta \,d\eta_1\left|q\left(\frac{\eta_1-\xi_1}{2}\right)\right|\frac{1}{4}(\eta_1+1)\|\widetilde q\|\\
\leqslant \frac{\|\widetilde q\|}{4^2}\int_0^\xi\,d\xi_1 \int_0^\eta \,d\eta_1 \left|q\left(\frac{\eta_1}{2}\right)\right| (\eta_1+1)\leqslant
\frac{\|\widetilde q\|^2}{4^2} \xi\frac{(\eta+2)^2}{2},
\end{eqnarray*}
where we use Lemma \ref{L1} to get the last inequality. Then the nest term is estimated by
\begin{eqnarray*}
|K^2Q(\xi,\eta)|\leqslant\frac{1}{4}\int_0^\xi\, d\xi_1 \int_\xi^\eta \,d\eta_1\left|q\left(\frac{\eta_1-\xi_1}{2}\right)\right|\frac{\|\widetilde q\|^2}{4^2} \xi_1\frac{(\eta_1+2)^2}{2}\\
\leqslant \frac{\|\widetilde q\|^2}{4^3}\int_0^\xi\,d\xi_1 \int_0^\eta \,d\eta_1 \left|q\left(\frac{\eta_1}{2}\right)\right| \xi_1\frac{(\eta_1+2)^2}{2}\leqslant
\frac{\|\widetilde q\|^3}{4^3} \frac{\xi^2}{2}\frac{(\eta+3)^3}{3!}.
\end{eqnarray*}
Using induction and Lemma \ref{L1} it is not hard to show that
\begin{eqnarray*}
|K^nQ(\xi,\eta)|\leqslant
\frac{\|\widetilde q\|^{n+1}}{4^{n+1}} \frac{\xi^n}{n!}\frac{(\eta+n+1)^{n+1}}{(n+1)!},
\end{eqnarray*}

Then for the general term of (\ref{gursa_sol_sum}) one can write an estimate using the inequality $(a+b)^n\leqslant 2^{n-1}\left(a^n+b^n\right)$:
\begin{eqnarray*}
|K^nQ(\xi,\eta)|\leqslant
\frac{\|\widetilde q\|^{n+1}}{4^{n+1}} \frac{\xi^n}{n!}\frac{(\eta+n+1)^{n+1}}{(n+1)!}\leqslant \frac{\|\widetilde q\|^{n+1}}{4^{n+1}} \frac{\xi^n}{n!}\frac{2^n\left(\eta^{n+1}+(n+1)^{n+1}\right)}{(n+1)!}\\
\leqslant \frac{\|\widetilde q\|^{n+1}}{2^{n+2}} \left(\frac{\xi^n}{n!}\frac{\eta^{n+1}}{(n+1)!}+\frac{\xi^n}{n!} \frac{(n+1)^{n+1}}{(n+1)!}\right).
\end{eqnarray*}
Using the above estimate and the Stirling inequality
\begin{equation*}
n!\geqslant \sqrt{2\pi}\left(\frac{n}{e}\right)^n,
\end{equation*}
and taking some $\varkappa>0$ we can proceed to the estimate of (\ref{gursa_sol_sum}):
\begin{eqnarray*}
|v(\xi,\eta)|\leqslant \sum_{n=0}^\infty \frac{\|\widetilde q\|^{n+1}}{2^{n+2}} \frac{\xi^n}{n!}\frac{\eta^{n+1}}{(n+1)!}+ \sum_{n=0}^\infty \frac{\|\widetilde q\|^{n+1}}{2^{n+2}}
\frac{\xi^n}{n!} \frac{(n+1)^{n+1}}{(n+1)!}\\
\leqslant  \sum_{n=0}^\infty \frac{\|\widetilde q\|\eta}{4(n+1)} \frac{\left(\frac{\xi \|\widetilde q\|\varkappa}{2} \right)^n}{n!}\frac{\left(\frac{\eta}{\varkappa}\right)^{n}}{n!}+
\sum_{n=0}^\infty \frac{\|\widetilde q\|e}{4 \sqrt{2\pi}}
\frac{\left(\frac{\xi e \|\widetilde q\|}{2}\right)^n}{n!} \\
\leqslant \frac{\|\widetilde q\|\eta}{4}e^{\frac{\xi \|\widetilde q\|\varkappa}{2}}e^{\frac{\eta}{\varkappa}}+\frac{\|\widetilde q\|e}{4 \sqrt{2\pi}}e^{\frac{\xi e \|\widetilde q\|}{2}}.
\end{eqnarray*}
Rewriting the above inequality for $w(x,t)$ we get that
\begin{equation}
\label{W_est2} |w(x,t)|\leqslant \frac{\|\widetilde q\|(x+t)}{4}
e^{\frac{(t-x)\|\widetilde
q\|\varkappa}{2}+\frac{t+x}{\varkappa}}+\frac{\|\widetilde q\|e}{4
\sqrt{2\pi}}e^{\frac{(t-x) e \|\widetilde q\|}{2}}.
\end{equation}
Then (\ref{W_est2}) implies that the Fourier transform
(\ref{Weyl_func}) of $u^\delta$  exists as soon as
$\operatorname{Im}k
>\operatorname{max}\left\{\frac{e\|\widetilde q\|}{2},\,
\frac{\varkappa\|\widetilde q\|}{2}\right\}$, and that
(\ref{Weyl_func}) is from $L_2(\mathbb{R}_+)$ if $\frac{\varkappa
\|\widetilde q\|}{2}>\frac{1}{\varkappa}$. Finally we arrive at
the following
\begin{proposition}
If $q\in l_\infty(l_1(\mathbb{R}_+))$ then formula
(\ref{Weyl_func}) gives the representation for the Weyl solution
of the Schr\"odinger operator $H$ in the region
 $\operatorname{Im}k >\operatorname{max}\left\{\frac{e\|\widetilde q\|}{2},\, \sqrt{\frac{|\widetilde q\|}{2}}\right\}$.
\end{proposition}

\end{document}